\documentclass[11pt,a4paper]{article}
\usepackage{fullpage}
\usepackage{authblk}

\usepackage{mathrsfs}
\usepackage{amsfonts}
\usepackage{amssymb}
\usepackage{bm}
\usepackage{multirow,bigstrut}
\usepackage{threeparttable}
\usepackage{enumerate}
\usepackage{tikz}
\usepackage{graphicx}
\usepackage{stfloats}
\usepackage{array}
\usepackage{amsmath}
\usepackage[colorlinks,
            linkcolor=red,
            anchorcolor=blue,
            citecolor=green
            ]{hyperref}

\newtheorem{Th}{Theorem}[section]

\newtheorem{Cor}{Corollary}[section]
\newtheorem{Pro}{Proposition}[section]

\newtheorem{Rem}{Remark}[section]
\newtheorem{Exa}{Example}[section]

\title{Quaternionic left eigenvalue problem: a matrix representation}
	\author{Wankai Liu\thanks{zjnulwk@163.com} }
	
	\author{Kit Ian Kou\thanks{kikou@umac.mo}}
	
	\affil{\normalsize{Department of Mathematics, Faculty of Science and Technology, University of Macau, Macao, China}}
	
	\date{}

\begin{document}
  \maketitle
\begin{abstract}
\normalsize

This paper presents an innovative set of tools developed to support a methodology to find the left eigenvalues of $m$ order quaternion square matrix. It is solving four real polynomial equations of order not greater than $4m-3$ in four variables. Some important properties of these eigenvalues are also investigated.
\end{abstract}


\begin{msc}
 42A15,
\end{msc}

\section{Introduction}\label{S1}
As usual, let $\mathbb{R}$, $\mathbb{H}$, $\mathbb{R}^{m \times n}$ and $\mathbb{H}^{m \times n}$ denote the sets of the real number, quaternion, $m\times n$ real matrix and $m \times n$ quaternion matrix. Quaternion is generally represented in the form
\[
q :=[q]_{0}+[q]_{1}\hbar+[q]_{2}\jmath+[q]_{3}\kappa \,\,\,\,
\]
with real coefficients $[q]_{0},[q]_{1},[q]_{2},[q]_{3}$ and  $\hbar^{2}=\jmath^{2}=\kappa^{2}=\hbar\jmath\kappa=-1$. Let $q^{\ast} :=[q]_{0}-[q]_{1}\hbar-[q]_{2}\jmath-[q]_{3}\kappa$ be the conjugate of $q$, and $\|q\| :=\sqrt{q^{\ast}q}=\sqrt{[q]_{0}^{2}+[q]_{1}^{2}+[q]_{2}^{2}+[q]_{3}^{2}}$, be the norm of $q$, and thus $q^{-1}=\frac{q^{\ast}}{q^{\ast}q}=\frac{q^{\ast}}{\|q\|^{2}}$. For any quaternion matrix $\mathbf{A} :=(a_{ij})\in \mathbb{H}^{m\times n}$, $1 \leq i\leq m$, $1 \leq j \leq n$, let $\mathbf{A}^{H} :=(a_{ij}^{\ast})^{T}$ be the conjugate transpose of $\mathbf{A}$ and Rank$(\mathbf{A})$ be the rank of $\mathbf{A}$. Two quaternions $a$ and $b$ are said to be similar if there exists a nonzero quaternion $\sigma$ such that $\sigma^{-1}a\sigma=b$, that is written as $a \sim b$.

Due to the multiplication of two quaternions is non-commutative, the left
$\lambda_{l}$ and the right eigenvalues $\lambda_{r}$ of the quaternion matrix $\mathbf{A}\in \mathbb{H}^{m \times m}$ need to be treated independently, that is
\begin{itemize}
  \item Right eigenvalue problem: $\mathbf{A}\mathbf{x}=\mathbf{x}\lambda_{r}$, \,\,for nonzero $\mathbf{x}\in \mathbb{H}^{m \times 1}$.

  \item Left eigenvalue problem: $\mathbf{A}\mathbf{x}=\lambda_{l} \mathbf{x}$, \,\,for nonzero $\mathbf{x}\in \mathbb{H}^{m \times 1}$.
\end{itemize}

The right eigenvalue problems are well established in \cite{brenner1951matrices,lee1948eigenvalues,zhang1997quaternions}, whereas the left eigenvalue problem is less solved.
The existence of left eigenvalues was proved by the topological method in \cite{wood1985quaternionic}. 
Huang \cite{huang2001left} explained how to compute all left eigenvalues of a $2\times 2$ matrix by the quadratic formulas for quaternions \cite{huang2002quadratic}. Macias \cite{macias2014topological} presented an incomplete classification of the left quaternion eigenvalue problems of $3\times3$ quaternion matrices by applying the characteristic map. It is still an open question for computing left eigenvalues for $n\times n$ quaternion matrices if $n\geq 3$.

Although So \cite{so2005quaternionic} showed that all left eigenvalues of $n\times n$ quaternion  matrix could be found by solving quaternion polynomials of degree not greater than $n$. Not only is it awkward to obtain the resulting quaternion polynomial, but there is not any known method to solve the resulting quaternion polynomial while some quaternion polynomial problems are studied \cite{kalantari2013algorithms,falcao2018evaluation}. In this work, we develop a novel method of computing left eigenvalues via solving four real polynomial equations of degree not greater than $4m-3$ with four variables. The contributions of this paper are summarised as follows.

\begin{enumerate}
  \item We propose the generalized characteristic polynomial of quaternion matrix. That is, the roots of the generalized characteristic polynomial are the left quaternionic eigenvalues.

 \item A condition of equivalence of the generalized characteristic polynomial is obtained, i.e. four real high order polynomial equations of four variables, the problem of computing the left eigenvalues, then, is to solve special polynomial equations.

  \item All the left eigenvalues of quaternion matrix are located in the in particular annulus connected with the right eigenvalues.
\end{enumerate}

The outline of the paper is as follows. In Section \ref{S2}, we describe the proposed method based on matrix representation. Section \ref{S3} gives some essential properties of left eigenvalue problems. Section \ref{S4} offers several illustrative examples. The last section \ref{S5} concludes.
\section{Calculation of left eigenvalues}\label{S2}
A total of 48 real matrix forms represent a quaternion \cite{farebrother2003matrix}. Let $\mathbf{E}$ be the $4 \times 4$ identity matrix and for $k=1,2,\cdots,48$, $\mathbf{H}_{k},\mathbf{J}_{k},\mathbf{K}_{k}$ be the $4\times4$ matrices with real entries. A quaternion $q$ can be written as various forms of real matrix
\begin{eqnarray}\label{48forms}
[q]_{0}\mathbf{E}+[q]_{1}\mathbf{H}_{k}+[q]_{2}\mathbf{J}_{k}+[q]_{3}\mathbf{K}_{k}\,\,,\,\,\,\, k=1,2,\cdots,48.
\end{eqnarray}
Quaternion addition and multiplication correspond to matrix addition and multiplication provided that the matrices $\mathbf{H}_{k},\mathbf{J}_{k},\mathbf{K}_{k}$ satisfying the "Hamiltonian conditions" ($\mathbf{H}_{k}\mathbf{H}_{k}=\mathbf{J}_{k}\mathbf{J}_{k}=\mathbf{K}_{k}\mathbf{K}_{k}=\mathbf{H}_{k}\mathbf{J}_{k}\mathbf{K}_{k}=-\mathbf{E}$).

First of all, we define two useful quaternion operators. For any quaternion $q\in\mathbb{H}$, we define the operator $\mathcal{Q}_{k}$ to maps a quaternion to the $4 \times 4$ real matrix. That is, $\mathcal{Q}_{k}: \mathbb{H}\rightarrow \mathbb{R}^{4\times 4}$, one defines
\begin{equation}
\mathcal{Q}_{k}(q)=[q]_{0}\mathbf{E}+[q]_{1}\mathbf{H}_{k}+[q]_{2}\mathbf{J}_{k}+[q]_{3}\mathbf{K}_{k}.
\,\,\,\,k=1,2,\cdots,48.\label{Q}
\end{equation}

In this representation (\ref{Q}), the conjugate of quaternion corresponds to the transpose of the matrix $\mathcal{Q}_{k}(q)$. The fourth power of the norm of a quaternion is the determinant of the corresponding matrix $\mathcal{Q}_{k}(q)$ (denote $\det\mathcal{Q}_{k}(q)$ ).

For any quaternion matrix $\mathbf{A}=(a_{ij})\in\mathbb{H}^{m\times n}$, we define the operator $\mathcal{P}_{k}$ to maps a $m \times n$ quaternion matrix to the $4m \times 4n$ real matrix. That is, $\mathcal{P}_{k}: \mathbb{H}^{m\times n}\rightarrow \mathbb{R}^{4m\times 4n}$, one defines
\begin{equation}
\mathcal{P}_{k}(\mathbf{A})=(\mathcal{Q}_{k}(a_{ij}))\,\,\,\,k=1,2,\cdots,48.
\label{P}
\end{equation}

Since quaternion addition and multiplication correspond to matrix addition and multiplication provided
that the matrices $\mathbf{E}$, $\mathbf{H}_{k}$, $\mathbf{J}_{k}$, $\mathbf{K}_{k}$. Some properties of operators $\mathcal{Q}_{k}$ and $\mathcal{P}_{k}$ are given as below.

\begin{Pro}\label{P2.1}
$\mathcal{Q}_{k}$ and $\mathcal{P}_{k}$ have the following properties
\begin{enumerate}
  \item $q_{1}+q_{2}=q_{3}$ if and only if $\mathcal{Q}_{k}(q_{1})+\mathcal{Q}_{k}(q_{2})=\mathcal{Q}_{k}(q_{3})$, $\forall\, q_{1},q_{2},q_{3}\in\mathbb{H}$.

  \item If $q_{1}q_{2}=q_{3}$ if and only if  $\mathcal{Q}_{k}(q_{1})\mathcal{Q}_{k}(q_{2})=\mathcal{Q}_{k}(q_{3})$, $\forall\, q_{1},q_{2},q_{3}\in\mathbb{H}$.

  \item $q\mathbf{A}=\mathbf{B}$ if and only if  $\mathcal{P}_{k}(q\mathbf{I})\mathcal{P}_{k}(\mathbf{A})=\mathcal{P}_{k}(\mathbf{B})$, $\forall\, q\in\mathbb{H}$, $\mathbf{I}$ is $m\times m$ identity matrix, $\forall$ \,$\mathbf{A}$, $\mathbf{B} \in \mathbb{H}^{m\times n}$.

  \item $\mathbf{A}q=\mathbf{B}$ if and only if $\mathcal{P}_{k}(\mathbf{A})\mathcal{P}_{k}(q\mathbf{I})=\mathcal{P}_{k}(\mathbf{B})$, $\forall\, q\in\mathbb{H}$, $\mathbf{I}$ is $n\times n$ identity matrix, $\forall$ \,$\mathbf{A}$, $\mathbf{B} \in \mathbb{H}^{m\times n}$.

  \item If $\mathbf{A}\mathbf{B}=\mathbf{C}$ if and only if $\mathcal{P}_{k}(\mathbf{A})\mathcal{P}_{k}(\mathbf{B})=\mathcal{P}_{k}(\mathbf{C})$, $\forall\, \mathbf{A}\in\mathbb{H}^{m\times r},\forall \, \mathbf{B}\in\mathbb{H}^{r\times n},\forall\, \mathbf{C}\in\mathbb{H}^{m\times n}$.
\end{enumerate}
\end{Pro}
\begin{Th}\label{2.1}
For any matrix $\mathbf{A}=(a_{ij})\in\mathbb{H}^{m\times n}$ and $m\geq n\geq1$, there exist nonnegative integer $s$, with $0\leq s\leq n$, such that $\mbox{Rank}(\mathcal{P}_{k}(\mathbf{A}))=4s$.
\end{Th}

{\bf Proof.}
We show that the statement $Rank(\mathcal{P}_{k}(\mathbf{A}))=4s$ by mathematical induction.

When $m=1$, the statement $Rank(\mathcal{P}_{k}(\mathbf{A}))=4s$ is straightforward for $s=0$ and $1$, since $\det{P}_{k}(\mathbf{A})=\det{Q}_{k}(a_{11})=\|a_{11}\|^{4}\geq0$.

Suppose that when $m=t$, statement holds, now we show that when $m=t+1$, the statement still holds. This can be done as follows. By Eq.(\ref{P}),
\begin{equation}\label{Pk}
\mathcal{P}_{k}(\mathbf{A}) =\left ( \begin{array}{cccc}  \mathcal{Q}_{k}(a_{11})&\mathcal{Q}_{k}(a_{12})&\cdots&\mathcal{Q}_{k}(a_{1n})
\\
\mathcal{Q}_{k}(a_{21})&\mathcal{Q}_{k}(a_{22})&\cdots&\mathcal{Q}_{k}(a_{2n})
\\
\vdots&\vdots&\cdots&\vdots
\\
\mathcal{Q}_{k}(a_{(t+1)1})&\mathcal{Q}_{k}(a_{(t+1)2})&\cdots&\mathcal{Q}_{k}(a_{(t+1)n})
\end{array}
\right ).
\end{equation}

If $\mathcal{P}_{k}(\mathbf{A})$ is a zero matrix, then the result is obvious. If it is nontrivial,
then there exist at least one entry $a_{i'j'}$ such that $\mathcal{Q}_{k}(a_{i'j'})\neq \mathbf{0}$. Without loss of generality (WLOG), let $\mathcal{Q}_{k}(a_{11})\neq \mathbf{0}$. Then there exist several block elementary transformations in row and column of the matrix such that
\begin{equation}\label{Qkrc}
\mathcal{P}_{k}(\mathbf{A})
\rightarrow\left ( \begin{array}{cccc}  \mathcal{Q}_{k}(a_{11})&\mathbf{0}&\cdots&\mathbf{0}
\\
\mathbf{0}&\mathcal{Q}_{k}(a^{\star}_{22})&\cdots&\mathcal{Q}_{k}(a^{\star}_{2n})
\\
\vdots&\vdots&\cdots&\vdots
\\
\mathbf{0}&\mathcal{Q}_{k}(a^{\star}_{(t+1)2})&\cdots&\mathcal{Q}_{k}(a^{\star}_{(t+1)n})
\end{array}
\right )
\end{equation}
and
\begin{equation}\label{Qkrcb}
\mathbf{B}=\left ( \begin{array}{ccc}

\mathcal{Q}_{k}(a^{\star}_{22})&\cdots&\mathcal{Q}_{k}(a^{\star}_{2n})
\\
\vdots&\cdots&\vdots
\\
\mathcal{Q}_{k}(a^{\star}_{(t+1)2})&\cdots&\mathcal{Q}_{k}(a^{\star}_{(t+1)n})
\end{array}
\right )
\end{equation}
where $\mathbf{0}=\mathcal{Q}_{k}(0)$ is the zero matrix. Using the induction hypothesis that for $m=t$ holds, then there exist nonnegative $s^{\star}$ ($0\leq s^{\star}\leq n-1$) such that $Rank(\mathbf{B})=4s^{\star}$.
According to (\ref{Qkrc}), we obtain $Rank(\mathcal{P}_{k}(\mathbf{A}))=4(s^{\star}+1)$ since $\mathcal{Q}_{k}(a_{11})\neq \mathbf{0}$. Thereby showing that indeed for $m=t+1$ holds since the  block elementary transformations do not change the rank of matrix. And the proof is complete.

\begin{Th}\label{2.2}
For any matrix $\mathcal{P}_{k}(\mathbf{A})=[\mathbf{a}_{1},\mathbf{a}_{2},\cdots,\mathbf{a}_{4n}]$ and $\mathcal{P}_{k}(\mathbf{B})=[\mathbf{b}_{1},\mathbf{b}_{2},\\ \mathbf{b}_{3},\mathbf{b}_{4}]$ corresponding  to $\mathbf{A}=(a_{ij})\in\mathbb{H}^{m\times n}$ and  $\mathbf{B}=(b_{ij})\in\mathbb{H}^{m\times 1}$ respectively, where $\mathbf{a}_{s}$ and $\mathbf{b}_{t}$ are the column vectors ($1\leq s\leq 4n$ and $1\leq t\leq 4$). If there was a column vector $\mathbf{b}_{t'}$ ($1\leq t'\leq 4$) can be lineally expressed by  $\mathbf{a}_{1},\mathbf{a}_{2},\cdots,\mathbf{a}_{4n}$, then the other three column vectors can also lineally expressed by  $\mathbf{a}_{1},\mathbf{a}_{2},\cdots,\mathbf{a}_{4n}$.
\end{Th}

{\bf Proof.} WLOG. let $\mathbf{b}_{t'}=\mathbf{b}_{1}$, there exist $\mathbf{x}=[x_{1},x_{2},\cdots,x_{4n}]^{T}\in\mathbb{R}^{4n\times 1}$ such that $b_{1}=\mathcal{P}_{k}(A)\mathbf{x}=\sum\limits_{i=1}^{4n}a_{i}x_{i}$.
Obviously, we can find $\mathcal{P}_{k}(\mathbf{c})=[\mathbf{c}_{1},\mathbf{c}_{2},\mathbf{c}_{3},\mathbf{c}_{4}]$ corresponding to $\mathbf{c}=(c_{ij})\in\mathbb{H}^{n\times 1}$ such that $\mathbf{c}_{1}=\mathbf{x}$. And we can get $\mathbf{b}_{2}=\mathcal{P}_{k}(\mathbf{A})\mathbf{c}_{2}$, $\mathbf{b}_{3}=\mathcal{P}_{k}(\mathbf{A})\mathbf{c}_{3}$ and $\mathbf{b}_{4}=\mathcal{P}_{k}(\mathbf{A})\mathbf{c}_{4}$.

Now we can obtain the following result according to the above Theorem \ref{2.1} and Theorem \ref{2.2}.
\begin{Th}\label{2.3} For any nonzero matrix $\mathbf{A}=(a_{ij})\in\mathbb{H}^{m\times m}$, if $Rank(\mathcal{P}_{k}(\mathbf{A}))=4n$, ($n=1,2,\dots,m$), then there exist $n$ order sub-matrix $\mathbf{B}\in\mathbb{H}^{n\times n}$ such that $Rank(\mathcal{P}_{k}(\mathbf{B}))=4n$.
\end{Th}

WLOG, we take
\[
\mathbf{H}_{1}=\left ( \begin{array}{cccc}  0  & -1 & 0 & 0
\\
1  & 0 & 0 & 0
\\
0  & 0 & 0 & -1
\\
0  & 0 & 1 & 0
\end{array}
\right )
\,\,\,\,
\mathbf{J}_{1}=\left ( \begin{array}{cccc}  0  & 0 & -1 & 0
\\
0  & 0 & 0 & 1
\\
1  & 0 & 0 & 0
\\
0  & -1 & 0 & 0
\end{array}
\right )
\,\,\,\,
\mathbf{K}_{1}=\left ( \begin{array}{cccc}  0  & 0 & 0 & -1
\\
0  & 0 & -1 & 0
\\
0  & -1 & 0 & 0
\\
-1 & 0 & 0 & 0
\end{array}
\right )
\]

then we obtain
\begin{equation}\label{Q11}
\mathcal{P}_{1}(\mathbf{A})=\left ( \begin{array}{cccc}  \mathcal{Q}_{1}(a_{11})&\mathcal{Q}_{1}(a_{12})&\cdots&\mathcal{Q}_{1}(a_{1m})
\\
\mathcal{Q}_{1}(a_{21})&\mathcal{Q}_{1}(a_{22})&\cdots&\mathcal{Q}_{1}(a_{2m})
\\
\vdots&\vdots&\cdots&\vdots
\\
\mathcal{Q}_{1}(a_{m1})&\mathcal{Q}_{1}(a_{m2})&\cdots&\mathcal{Q}_{1}(a_{mm})
\end{array}
\right )
\end{equation}
\begin{Th}\label{2.4} For any matrix $\mathbf{A}=(a_{ij})\in\mathbb{H}^{m\times m}$, $m>1$, then the following statements are equivalent:

\begin{itemize}
\item[(1).] $\det\mathcal{P}_{1}(\mathbf{A})=0$.

\item[(2).] The determinants of all $4m-3$ order sub-matrix of $\mathcal{P}_{1}(\mathbf{A})$ are zero.

\item[(3).] Let $\mathcal{S}=\{\mathbf{B}:\mathbf{B}\,\,is\,\,m-1\,\,order\,\,submatrix\,\,of\,\,\mathbf{A}\}$ and  $\max\{Rank(\mathcal{P}_{1}(\mathbf{B})):\mathbf{B}\in \mathcal{S}\}=4n$ ($0\leq n\leq m-1$). There exist the determinants of sixteen $4m-3$ order sub-matrix  (denoted as $\mathbf{C}_{t}$, $t=1,2,\cdots,16$) of $\mathcal{P}_{1}(\mathbf{A})$ are zero, where $\mathcal{P}_{1}(\mathbf{\mathbf{B}}')$ is a sub-matrix of $\mathbf{C}_{t}$ and $Rank(\mathcal{P}_{1}(\mathbf{B}'))=4n$ ($\mathbf{B}'\in \mathcal{S}$).
\end{itemize}
\end{Th}

{\bf Proof.} $1\Leftrightarrow 2$: this can be obtained according to the above Theorem \ref{2.1}.

$1\Rightarrow3$: this can be seen from $1\Leftrightarrow2$.

$3\Rightarrow1$: if $n<m-1$, this can be obtained by the above Theorem \ref{2.2}. If $n=m-1$, WLOG, let
\[
\mathbf{B}'=\left(\begin{array}{ccc}
\mathcal{Q}_{1}(a_{22})&\cdots&\mathcal{Q}_{1}(a_{2m})
\\
\vdots&\cdots&\vdots
\\
\mathcal{Q}_{1}(a_{m2})&\cdots&\mathcal{Q}_{1}(a_{mm})
\end{array}
\right)
\]
and
\begin{equation}\label{ct}
\begin{array}{cc}
\mathbf{C}_{1}=(\mathcal{P}_{1}(\mathbf{A}))(1,5:4m;1,5:4m)&\mathbf{C}_{2}=(\mathcal{P}_{1}(\mathbf{A}))(2,5:4m;2,5:4m)
\\
\mathbf{C}_{3}=(\mathcal{P}_{1}(\mathbf{A}))(3,5:4m;3,5:4m)&\mathbf{C}_{4}=(\mathcal{P}_{1}(\mathbf{A}))(4,5:4m;4,5:4m)
\\
\mathbf{C}_{5}=(\mathcal{P}_{1}(\mathbf{A}))(1,5:4m;2,5:4m)&\mathbf{C}_{6}=(\mathcal{P}_{1}(\mathbf{A}))(2,5:4m;1,5:4m)
\\
\mathbf{C}_{7}=(\mathcal{P}_{1}(\mathbf{A}))(3,5:4m;4,5:4m)&\mathbf{C}_{8}=(\mathcal{P}_{1}(\mathbf{A}))(4,5:4m;3,5:4m)
\\
\mathbf{C}_{9}=(\mathcal{P}_{1}(\mathbf{A}))(1,5:4m;3,5:4m)&\mathbf{C}_{10}=(\mathcal{P}_{1}(\mathbf{A}))(2,5:4m;4,5:4m)
\\
\mathbf{C}_{11}=(\mathcal{P}_{1}(\mathbf{A}))(3,5:4m;1,5:4m)&\mathbf{C}_{12}=(\mathcal{P}_{1}(\mathbf{A}))(4,5:4m;2,5:4m)
\\
\mathbf{C}_{13}=(\mathcal{P}_{1}(\mathbf{A}))(1,5:4m;4,5:4m)&\mathbf{C}_{14}=(\mathcal{P}_{1}(\mathbf{A}))(2,5:4m;3,5:4m)
\\
\mathbf{C}_{15}=(\mathcal{P}_{1}(\mathbf{A}))(3,5:4m;2,5:4m)&\mathbf{C}_{16}=(\mathcal{P}_{1}(\mathbf{A}))(4,5:4m;2,5:4m)
\end{array}
\end{equation}

where $(\mathcal{P}_{1}(\mathbf{A}))(r,5:4m;c,5:4m)$ is obtained by taking row $r$ and row from $5$ to $4n$ and column $c$ and column from $5$ to $4n$ of $\mathcal{P}_{1}(\mathbf{A})$.

If $\det\mathbf{C}_{t}=0$, then $\mathbf{C}_{4}(:,1)$ can be lineally expressed by the columns of $\mathbf{C}_{4}(:,2:4n+1)$ since $Rank(\mathcal{P}_{k}(\mathbf{B}'))=4n$, i.e. there exist column vector $\mathbf{x}_{1}\in\mathbb{H}^{4n\times 1}$ such that $\mathbf{C}_{4}(:,1)=\mathbf{C}_{4}(:,2:4n+1)\mathbf{x}_{1}$. Similarly, there exist column vector $\mathbf{x}_{2}$ $\mathbf{x}_{3}$ $\mathbf{x}_{4}$ such that $\mathbf{C}_{7}(:,1)=\mathbf{C}_{7}(:,2:4n+1)\mathbf{x}_{2}$, $\mathbf{C}_{10}(:,1)=\mathbf{C}_{10}(:,2:4n+1)\mathbf{x}_{3}$ and $\mathbf{C}_{13}(:,1)=\mathbf{C}_{13}(:,2:4n+1)\mathbf{x}_{4}$. And because $Rank(\mathcal{Q}_{k}(B'))=4n$, we get $\mathbf{x}_{1}=\mathbf{x}_{2}=\mathbf{x}_{3}=\mathbf{x}_{4}$. Therefore, $(\mathcal{P}_{1}(\mathbf{A}))(:,4)$ can be lineally expressed by the columns of $(\mathcal{P}_{1}(\mathbf{A}))(:,5:4m)$, i.e. $\det\mathcal{P}_{k}(\mathbf{A})=0$.

\begin{Th}\label{2.5} For any quaternion $q$, then

\begin{itemize}
\item[(1).]$\mathbf{P}(3,4)\mathbf{P}(1,2)\mathbf{P}(4(-1))\mathbf{P}(1(-1))\mathcal{P}_{1}(q)\mathbf{P}(1(-1))\mathbf{P}(4(-1))\mathbf{P}(1,2)\mathbf{P}(3,4)=\mathcal{P}_{1}(q)$
\item[(2).]$\mathbf{P}(2,4)\mathbf{P}(1,3)\mathbf{P}(4(-1))\mathbf{P}(3(-1))\mathcal{P}_{1}(q)\mathbf{P}(3(-1))\mathbf{P}(4(-1))\mathbf{P}(1,3)\mathbf{P}(2,4)=\mathcal{P}_{1}(q)$
\item[(3).]$\mathbf{P}(2,3)\mathbf{P}(1,4)\mathbf{P}(4(-1))\mathbf{P}(2(-1))\mathcal{P}_{1}(q)\mathbf{P}(2(-1))\mathbf{P}(4(-1))\mathbf{P}(1,4)\mathbf{P}(2,3)=\mathcal{P}_{1}(q)$
\end{itemize}
\end{Th}
where the elementary matrix $\mathbf{P}(i,j)$ is obtained by swapping row $i$ and row $j$ of the $4\times 4$ identity matrix and $\mathbf{P}(i(c))$ is a $4\times 4$ diagonal matrix, with diagonal entries $1$ everywhere except in the $ith$  position, where it is $c$.

Using the Theorem \ref{2.5}, it is not hard to see the following

\begin{Th}\label{2.6}
For any matrix $\mathbf{A}=(a_{ij})\in\mathbb{H}^{m\times m}$, $m>1$, then
\[
\det\mathbf{C}_{1}=\det\mathbf{C}_{2}=\det\mathbf{C}_{3}=\det\mathbf{C}_{4},\,\, -\det\mathbf{C}_{5}=\det\mathbf{C}_{6}=-\det\mathbf{C}_{7}=\det\mathbf{C}_{8},\,\,
\]
\[
-\det\mathbf{C}_{9}=\det\mathbf{C}_{10}=\det\mathbf{C}_{11}=-\det\mathbf{C}_{12},\,\,
-\det\mathbf{C}_{13}=-\det\mathbf{C}_{14}=\det\mathbf{C}_{15}=\det\mathbf{C}_{16}.
\]
\end{Th}
where $\mathbf{C}_{t}$ is defined in a similar way described in the Theorem \ref{2.4}.

Using Theorems \ref{2.4} and \ref{2.6}, it is not hard to see the following.

\begin{Th}\label{2.7}
For any matrix $\mathbf{A}=(a_{ij})\in\mathbb{H}^{m\times m}$, $m>1$, then the following statements are equivalent:
\begin{itemize}
\item[(1).] $\det\mathcal{P}_{k}(\mathbf{A})=0$
\item[(2).] $\exists\, t_{1}\in\{1,2,3,4\}$, $\exists\, t_{2}\in\{5,6,7,8\}$, $\exists\, t_{3}\in\{9,10,11,12\}$, and $\exists\, t_{4}\in\{13,14,15,16\}$, such that $det\mathbf{C}_{t_{1}}=det\mathbf{C}_{t_{2}}=det\mathbf{C}_{t_{3}}=det\mathbf{C}_{t_{4}}=0$
\item[(3).] $\forall\, t_{1}\in\{1,2,3,4\}$, $\forall\, t_{2}\in\{5,6,7,8\}$, $\forall\, t_{3}\in\{9,10,11,12\}$, and $\forall\, t_{4}\in\{13,14,15,16\}$, such that $det\mathbf{C}_{t_{1}}=det\mathbf{C}_{t_{2}}=det\mathbf{C}_{t_{3}}=det\mathbf{C}_{t_{4}}=0$
\end{itemize}
%
\end{Th}
where $\mathbf{C}_{t}$ is defined in a similar way described in the Theorem \ref{2.4}.

Let $\lambda_{l}\in \mathbb{H}$ be a left eigenvalue of $\mathbf{A}=(a_{ij})\in\mathbb{H}^{m\times m}$ with eigenvector $\mathbf{v}\in \mathbb{H}^{m}$, by the property $(5)$ of Proposition \ref{P2.1}, get
\begin{equation}
\mathbf{A}\mathbf{v}=\lambda_{l}\mathbf{v}\Leftrightarrow(\mathbf{A}-\lambda_{l}\mathbf{I})\mathbf{v}=0\Leftrightarrow\mathcal{P}_{1}(\mathbf{A}-\lambda_{l}\mathbf{I})\mathcal{P}_{1}(\mathbf{v})=\mathcal{P}_{1}(0).
\label{a3}
\end{equation}
Since $\mathbf{v}\in \mathbb{H}^{m}$ is nonzero, then the following homogenous linear equations
\begin{equation}
\mathcal{P}_{1}(\mathbf{A}-\lambda_{l}\mathbf{I})\mathbf{y}=0, \,\,\,\,\mathbf{y}\in \mathbb{R}^{4m\times 1}.
\end{equation}
have nonzero solution, i.e. $\det\mathcal{P}_{1}(\mathbf{A}-\lambda_{l}\mathbf{I})=0$. And $\det\mathcal{P}_{1}(\mathbf{A}-\lambda_{l}\mathbf{I})$ is called the generalized characteristic polynomial of $\mathbf{A}$. By the Theorem \ref{2.7}, We can solve and analyze the left eigenvalues from the specified four real $4m-3$ order polynomial equations in four variables. (Theorem \ref{2.7} $(2)$ ).

\section{Some properties}\label{S3}

\begin{Th}\label{3.1}
Let $\lambda_{l}\in \mathbb{H}$ is a left ($\lambda_{r}$ right) eigenvalue of matrix $\mathbf{A}=(a_{ij})\in\mathbb{H}^{m\times m}$, for any matrix $\mathbf{B}=(b_{ij})\in\mathbb{H}^{m\times m}$, if there exist  $\widetilde{k}$ such that $\mathcal{P}_{k}(\mathbf{A})=\mathcal{P}_{\widetilde{k}}(\mathbf{B})$, then $\gamma_{l}$ is a left ($\gamma_{r}$ right) eigenvalue of matrix $\mathbf{B}$, if $\mathcal{Q}_{\widetilde{k}}(\gamma_{l})=\mathcal{Q}_{k}(\lambda_{l})$ ($\mathcal{Q}_{\widetilde{k}}(\gamma_{r})=\mathcal{Q}_{k}(\lambda_{r})$).
\end{Th}
 {\bf Proof.} Let $\lambda_{l}$ be a left eigenvalue of $\mathbf{A}$ with eigenvector $\mathbf{v}$, then
 \[
 (\mathcal{P}_{\widetilde{k}}(\mathbf{B})-\mathcal{P}_{\widetilde{k}}(\gamma_{l} \mathbf{I}))\mathcal{P}(\mathbf{v})=(\mathcal{P}_{k}(\mathbf{A})-\mathcal{P}_{k}(\lambda_{l} \mathbf{I}))\mathcal{P}(\mathbf{v})=0
 \]
 by the Eq.\ref{a3}, the conclusion then follows immediately.

If $\lambda_{r}$ be a right eigenvalue of $\mathbf{A}$ with eigenvector $\mathbf{v}$, then we can take nonzero vector $\mathbf{u}$ such that $\mathcal{P}_{\widetilde{k}}(\mathbf{u})=\mathcal{P}_{k}(\mathbf{v})$. By the Eq.\ref{a3}, get
\[
\mathcal{P}_{k}(A)\mathcal{P}_{k}(\mathbf{v})=\mathcal{P}_{k}(\mathbf{v})\mathcal{P}_{k}(\lambda_{r})\Leftrightarrow \mathcal{P}_{\widetilde{k}}(\mathbf{B})\mathcal{P}_{\widetilde{k}}(\mathbf{u})=\mathcal{P}_{\widetilde{k}}(\mathbf{u})\mathcal{P}_{\widetilde{k}}(\gamma_{r})\Leftrightarrow \mathbf{B}\mathbf{u}=\mathbf{u}\gamma_{r}
\]
The proof is complete.

\begin{Cor}\label{3.111}
If $\mathcal{P}_{k}(\mathbf{A})=\mathcal{P}_{\widetilde{k}}(\mathbf{B})$, then $\mathbf{A}$ and $\mathbf{B}$ are similar.
\end{Cor}

\begin{Th}\label{3.2}
Let $\lambda_{l}\in \mathbb{H}$ is a left eigenvalue of matrix $\mathbf{A}=(a_{ij})\in\mathbb{H}^{m\times m}$ with eigenvector $\mathbf{v}$, for any two nonzero quaternions $a$ and $b$, then $a \lambda_{l} b$ is a left eigenvalue of matrix $a\mathbf{A}b$.
\end{Th}
{\bf Proof.} Since $a,b$ are nonzero, then
\[
(\mathbf{A}-\lambda_{l}\mathbf{I})\mathbf{v}=0\Leftrightarrow a(\mathbf{A}-\lambda_{l}\mathbf{I})b(b^{-1}\mathbf{v})=0\Leftrightarrow (a\mathbf{A}b-a\lambda_{l}b\mathbf{I})(b^{-1}\mathbf{v})=0
\]

\begin{Th}\label{3.3}
Let $\lambda_{l}\in \mathbb{H}$ is a left eigenvalue of matrix $\mathbf{A}=(a_{ij})\in\mathbb{H}^{m\times m}$, then there exist nonnegative real numbers $\mathfrak{m}\geq 0$ and $\mathfrak{M}\geq 0$ such that $\mathfrak{m}\leq\|\lambda_{l}\|\leq \mathfrak{M}$.
\end{Th}
{\bf Proof.} Let $\lambda_{l}$ be a left eigenvalue of $\mathbf{A}$ with unit eigenvector $\mathbf{v}$, then $\|\lambda_{l}\|^{2}$$=\mathbf{v}^{H}\|\lambda_{l}\|^{2}\mathbf{v}$$=\mathbf{v}^{H}\lambda_{l}^{H}\lambda_{l}\mathbf{v}$$=\mathbf{v}^{H}\mathbf{A}^{H}\mathbf{Av}$. Using Singular-value decomposition theorem \cite{zhang1997quaternions}, there exist unitary quaternionic matrices $\mathbf{V}$ such that $\mathbf{A}$
\[
\mathbf{v}^{H}\mathbf{A}^{H}\mathbf{A}\mathbf{v}=\mathbf{v}^{H}\mathbf{V}^{H}\left ( \begin{array}{ll}   \mathbf{D}^{2}  & 0
 \\
0  & 0
\end{array}
\right )\mathbf{V}\mathbf{v}
\]
where $\mathbf{D}=diag\{d_{1},d_{2},\cdots,d_{r}\}$ and the $d$'s are the positive singular values of $\mathbf{A}$. Let $\mathbf{y}=\mathbf{V}\mathbf{v}=[y_{1},y_{2},\cdots,y_{n}]^{T}$, since $\mathbf{V}$ is a unitary matrix and $\mathbf{v}$ is unit eigenvector, then $\mathbf{y}$ is unit vector and
\[
\|\lambda\|^{2}=\mathbf{y}^{H}\left ( \begin{array}{ll}   D^{2}  & 0
 \\
0  & 0
\end{array}
\right )\mathbf{y}=\sum\limits_{l=1}^{r}d^{2}_{i}\|y_{i}\|^{2}
\]
if $r=n$, then
\[
d^{2}_{min}=d^{2}_{min}\sum\limits_{l=1}^{r}\|y_{i}\|^{2}\leq\|\lambda\|^{2}\leq d^{2}_{max}\sum\limits_{l=1}^{r}\|y_{i}\|^{2}=d^{2}_{max}
\]
otherwise, may have
\[
0\leq\|\lambda\|^{2}\leq d^{2}_{max}
\]
where $d_{min}=min\{d_{1},d_{2},\cdots,d_{r}\}$ and $d_{max}=max\{d_{1},d_{2},\cdots,d_{r}\}$. The conclusion then follows immediately.

It was given that the left spectrum is compact \cite{macias2014topological} and all the left eigenvalues of quaternion matrix are located in the union of $n$ Gersgorin balls \cite{zhang2007gervsgorin}. In fact, all the left eigenvalues of quaternion matrix are located in the in particular annulus connected with the singular values.

\begin{Rem}
the norm of left eigenvalues $\lambda_{l}$ is dominated by the norm of the right eigenvalues $\lambda_{r}$, i.e.
\[
\alpha\leq\|\lambda_{l}\|\leq\beta
\]
where $\alpha=min\{\|\lambda_{r}\|: \mathbf{Av}=\mathbf{v}\lambda_{r},\,\,\mathbf{v}\neq0\}$ and $\beta=max\{\|\lambda_{r}\|: \mathbf{Av}=\mathbf{v}\lambda_{r},\,\,\mathbf{v}\neq0\}$.
\end{Rem}

\section{Examples}\label{S4}

\begin{Exa}\label{4.1}
Let
\[
\mathbf{A}=
\left[
\begin{array}{cc}
1&\kappa
\\
\kappa&1
\end{array}
\right].
\]
\end{Exa}
Then
\begin{align}
\mathcal{P}_{1}(\mathbf{A}-\lambda_{l} \mathbf{I})=
\begin{bmatrix}\begin{smallmatrix}
1-[\lambda_{l}]_{0}&[\lambda_{l}]_{1}&[\lambda_{l}]_{2}&[\lambda_{l}]_{3}&0&0&0&-1
\\
-[\lambda_{l}]_{1}&1-[\lambda_{l}]_{0}&[\lambda_{l}]_{3}&-[\lambda_{l}]_{2}&0&0&-1&0
\\
-[\lambda_{l}]_{2}&-[\lambda_{l}]_{3}&1-[\lambda_{l}]_{0}&[\lambda_{l}]_{1}&0&1&0&0
\\
-[\lambda_{l}]_{3}&[\lambda_{l}]_{2}&-[\lambda_{l}]_{1}&1-[\lambda_{l}]_{0}&1&0&0&0
\\
0&0&0&-1&1-[\lambda_{l}]_{0}&[\lambda_{l}]_{1}&[\lambda_{l}]_{2}&[\lambda_{l}]_{3}
\\
0&0&-1&0&-[\lambda_{l}]_{1}&1-[\lambda_{l}]_{0}&[\lambda_{l}]_{3}&-[\lambda_{l}]_{2}
\\
0&1&0&0&-[\lambda_{l}]_{2}&-[\lambda_{l}]_{3}&1-[\lambda_{l}]_{0}&[\lambda_{l}]_{1}
\\
1&0&0&0&-[\lambda_{l}]_{3}&[\lambda_{l}]_{2}&-[\lambda_{l}]_{1}&1-[\lambda_{l}]_{0}
\end{smallmatrix}\end{bmatrix}
\end{align}
Let $\mathbf{B}=\mathcal{P}_{1}(\mathbf{A}-\lambda_{l} \mathbf{I})$, From the Theorem \ref{2.7}, we can take
\begin{equation}\label{4.1.1}
\det\mathbf{B}(1:5,4:8)=-[\lambda_{l}]_{0}^{2}+2[\lambda_{l}]_{0}-[\lambda_{l}]_{1}^{2}-[\lambda_{l}]_{2}^{2}+[\lambda_{l}]_{3}^{2}-2=0.
\end{equation}
\begin{equation}\label{4.1.2}
\det\mathbf{B}(1:4,6,4:8)=-2[\lambda_{l}]_{2}[\lambda_{l}]_{3}=0.
\end{equation}
\begin{equation}\label{4.1.3}
\det\mathbf{B}(1:4,7,4:8)=2[\lambda_{l}]_{1}[\lambda_{l}]_{3}=0.
\end{equation}
\begin{equation}\label{4.1.4}
\det\mathbf{B}(1:4,8,4:8)=-[\lambda_{l}]_{3}([\lambda_{l}]_{0}-1)=0.
\end{equation}
If $[\lambda_{l}]_{3}=0$, by the Eq.\ref{4.1.1}, \ref{4.1.2}, \ref{4.1.3} and \ref{4.1.4}, then $\lambda_{l}$ need satisfy
\[
\det\mathbf{B}(1:5,4:8)=-([\lambda_{l}]_{0}-1)^{2}-[\lambda_{l}]_{1}^{2}-[\lambda_{l}]_{2}^{2}-1=0.
\]
Obviously, a real solution of this equation does not exist.
\\
If $[\lambda_{l}]_{3}\neq0$, by the Eq.\ref{4.1.1}, \ref{4.1.2}, \ref{4.1.3} and \ref{4.1.4}, then
\[
[\lambda_{l}]_{3}^{2}-1=0,\,\,\,\,
[\lambda_{l}]_{2}=0,\,\,\,\,
[\lambda_{l}]_{1}=0,\,\,\,\,
[\lambda_{l}]_{0}=1.
\]
i.e.
\[
[\lambda_{l}]_{3}=\pm1,\,\,\,\,
[\lambda_{l}]_{2}=0,\,\,\,\,
[\lambda_{l}]_{1}=0,\,\,\,\,
[\lambda_{l}]_{0}=1.
\]

In conclusion, $\lambda_{l}=1\pm \kappa$.

\begin{Exa}\label{4.2}
Let
\[
\mathbf{A}=
\left[
\begin{array}{cc}
1&-\jmath
\\
\jmath&1
\end{array}
\right].
\]
\end{Exa}
Then
\begin{align}
\mathcal{P}_{1}(\mathbf{A}-\lambda_{l} \mathbf{I})=
\begin{bmatrix}\begin{smallmatrix}
1-[\lambda_{l}]_{0}&[\lambda_{l}]_{1}&[\lambda_{l}]_{2}&[\lambda_{l}]_{3}&0&0&1&0
\\
-[\lambda_{l}]_{1}&1-[\lambda_{l}]_{0}&[\lambda_{l}]_{3}&-[\lambda_{l}]_{2}&0&0&0&-1
\\
-[\lambda_{l}]_{2}&-[\lambda_{l}]_{3}&1-[\lambda_{l}]_{0}&[\lambda_{l}]_{1}&-1&0&0&0
\\
-[\lambda_{l}]_{3}&[\lambda_{l}]_{2}&-[\lambda_{l}]_{1}&1-[\lambda_{l}]_{0}&0&1&0&0
\\
0&0&-1&0&1-[\lambda_{l}]_{0}&[\lambda_{l}]_{1}&[\lambda_{l}]_{2}&[\lambda_{l}]_{3}
\\
0&0&0&1&-[\lambda_{l}]_{1}&1-[\lambda_{l}]_{0}&[\lambda_{l}]_{3}&-[\lambda_{l}]_{2}
\\
1&0&0&0&-[\lambda_{l}]_{2}&-[\lambda_{l}]_{3}&1-[\lambda_{l}]_{0}&[\lambda_{l}]_{1}
\\
0&-1&0&0&-[\lambda_{l}]_{3}&[\lambda_{l}]_{2}&-[\lambda_{l}]_{1}&1-[\lambda_{l}]_{0}
\end{smallmatrix}\end{bmatrix}
\end{align}
Let $\mathbf{B}=\mathcal{P}_{1}(\mathbf{A}-\lambda_{l} \mathbf{I})$, From the Theorem \ref{2.7}, we can take
\begin{equation}\label{4.2.1}
\det\mathbf{B}(1:5,4:8)=-2[\lambda_{l}]_{2}[\lambda_{l}]_{3}=0
\end{equation}
\begin{equation}\label{4.2.2}
\det\mathbf{B}(1:4,6,4:8)=-[\lambda_{l}]_{0}^{2}+2[\lambda_{l}]_{0}-[\lambda_{l}]_{1}^{2}+[\lambda_{l}]_{2}^{2}-[\lambda_{l}]_{3}^{2}=0
\end{equation}
\begin{equation}\label{4.2.3}
\det\mathbf{B}(1:4,7,4:8)=-2[\lambda_{l}]_{1}[\lambda_{l}]_{2}=0
\end{equation}
\begin{equation}\label{4.2.4}
\det\mathbf{B}(1:4,8,4:8)=[\lambda_{l}]_{2}([\lambda_{l}]_{0}-1)=0
\end{equation}
If $[\lambda_{l}]_{2}=0$, by the Eq.\ref{4.2.1}, \ref{4.2.2}, \ref{4.2.3} and \ref{4.2.4}, then $\lambda_{l}$ need satisfy
\[
-([\lambda_{l}]_{0}-1)^{2}-[\lambda_{l}]_{1}^{2}-[\lambda_{l}]_{3}^{2}+1=0
\]
i.e.
\[
([\lambda_{l}]_{0}-1)^{2}+[\lambda_{l}]_{1}^{2}+[\lambda_{l}]_{3}^{2}=1
\]
thus
\[
\lambda_{l}=[\lambda_{l}]_{0}+[\lambda_{l}]_{1}\hbar+[\lambda_{l}]_{3}\kappa,\,\,\,where\,\,\,([\lambda_{l}]_{0}-1)^{2}+[\lambda_{l}]_{1}^{2}+[\lambda_{l}]_{3}^{2}=1
\]
\\
If $[\lambda_{l}]_{2}\neq0$, by the Eq.\ref{4.2.1}, \ref{4.2.2}, \ref{4.2.3} and \ref{4.2.4}, then
\[
[\lambda_{l}]_{3}=0,\,\,\,\,
[\lambda_{l}]_{2}^{2}+1=0,\,\,\,\,
[\lambda_{l}]_{1}=0,\,\,\,\,
[\lambda_{l}]_{0}=1
\]
Obviously, a real solution of this equation does not exist.

In conclusion, $\lambda_{l}=[\lambda_{l}]_{0}+[\lambda_{l}]_{1}\hbar+[\lambda_{l}]_{3}\kappa$, where $([\lambda_{l}]_{0}-1)^{2}+[\lambda_{l}]_{1}^{2}+[\lambda_{l}]_{3}^{2}=1$.

In fact, the left eigenvalues of Example \ref{4.1} and \ref{4.2} can also be computed using the formula \cite{huang2001left} and the same result is obtained.


\begin{Exa}\label{4.3}
Let
\[
\mathbf{A}=
\left[
\begin{array}{ccc}
\kappa&0&1
\\
0&\kappa&0
\\
1&0&\kappa
\end{array}
\right].
\]
\end{Exa}
Then
\[
\mathcal{P}_{1}(\mathbf{A}-\lambda_{l} \mathbf{I})=
\]

\begin{align}
\begin{bmatrix}\begin{smallmatrix}
-[\lambda_{l}]_{0}&[\lambda_{l}]_{1}&[\lambda_{l}]_{2}&-1+[\lambda_{l}]_{3}&0&0&0&0&1&0&0&0
\\
-[\lambda_{l}]_{1}&-[\lambda_{l}]_{0}&-1+[\lambda_{l}]_{3}&-[\lambda_{l}]_{2}&0&0&0&0&0&1&0&0
\\
-[\lambda_{l}]_{2}&1-[\lambda_{l}]_{3}&-[\lambda_{l}]_{0}&[\lambda_{l}]_{1}&0&0&0&0&0&0&1&0
\\
1-[\lambda_{l}]_{3}&[\lambda_{l}]_{2}&-[\lambda_{l}]_{1}&-[\lambda_{l}]_{0}&0&0&0&0&0&0&0&1
\\
0&0&0&0&-[\lambda_{l}]_{0}&[\lambda_{l}]_{1}&[\lambda_{l}]_{2}&-1+[\lambda_{l}]_{3}&0&0&0&0
\\
0&0&0&0&-[\lambda_{l}]_{1}&-[\lambda_{l}]_{0}&-1+[\lambda_{l}]_{3}&-[\lambda_{l}]_{2}&0&0&0&0
\\
0&0&0&0&-[\lambda_{l}]_{2}&1-[\lambda_{l}]_{3}&-[\lambda_{l}]_{0}&[\lambda_{l}]_{1}&0&0&0&0
\\
0&0&0&0&1-[\lambda_{l}]_{3}&[\lambda_{l}]_{2}&-[\lambda_{l}]_{1}&-[\lambda_{l}]_{0}&0&0&0&0
\\
1&0&0&0&0&0&0&0&-[\lambda_{l}]_{0}&[\lambda_{l}]_{1}&[\lambda_{l}]_{2}&-1+[\lambda_{l}]_{3}
\\
0&1&0&0&0&0&0&0&-[\lambda_{l}]_{1}&-[\lambda_{l}]_{0}&-1+[\lambda_{l}]_{3}&-[\lambda_{l}]_{2}
\\
0&0&1&0&0&0&0&0&-[\lambda_{l}]_{2}&1-[\lambda_{l}]_{3}&-[\lambda_{l}]_{0}&[\lambda_{l}]_{1}
\\
0&0&0&1&0&0&0&0&1-[\lambda_{l}]_{3}&[\lambda_{l}]_{2}&-[\lambda_{l}]_{1}&-[\lambda_{l}]_{0}
\end{smallmatrix}\end{bmatrix}
\end{align}

Let $\mathbf{B}=\mathcal{P}_{1}(\mathbf{A}-\lambda_{l} \mathbf{I})$, From the Theorem \ref{2.7}, let $\lambda_{l}\neq \kappa$, we can take
\begin{equation}\label{4.3.1}
\det\mathbf{B}(1:9,4:12)=2[\lambda_{l}]_{0}([\lambda_{l}]_{3}-1)\alpha=0
\end{equation}
\begin{equation}\label{4.3.2}
\det\mathbf{B}(1:8,10,4:12)=-2[\lambda_{l}]_{0}[\lambda_{l}]_{2}\alpha=0
\end{equation}
\begin{equation}\label{4.3.3}
\det\mathbf{B}(1:8,11,4:12)=2[\lambda_{l}]_{1}[\lambda_{l}]_{1}\alpha=0
\end{equation}
\begin{equation}\label{4.3.4}
\det\mathbf{B}(1:8,12,4:12)=(-[\lambda_{l}]_{0}^{2}+[\lambda_{l}]_{1}^{2}+[\lambda_{l}]_{2}^{2}+[\lambda_{l}]_{3}^{2}-2[\lambda_{l}]_{3}+2)\alpha=0
\end{equation}
and
\[
\alpha=([\lambda_{l}]_{0})^{2}+([\lambda_{l}]_{1})^{2}+[\lambda_{l}]_{2}^{2}+([\lambda_{l}]_{3}-1)^{2}.
\]
Since $\lambda_{l}\neq \kappa$, thus
\[
\alpha\neq 0
\]
If $[\lambda_{l}]_{0}\neq0$, by the Eq.\ref{4.3.1}, \ref{4.3.2}, \ref{4.3.3} and \ref{4.3.4}, then
\[
[\lambda_{l}]_{1}=0,\,\,\,\,[\lambda_{l}]_{2}=0,\,\,\,\,[\lambda_{l}]_{3}=1\,\,\,\,and\,\,\,\,[\lambda_{l}]_{0}^{2}-1=0
\]
i.e.
\[
[\lambda_{l}]_{1}=0,\,\,\,\,[\lambda_{l}]_{2}=0,\,\,\,\,[\lambda_{l}]_{3}=1\,\,\,\,and\,\,\,\,[\lambda_{l}]_{0}=\pm1
\]
thus
\[
\lambda_{l}=\pm1+\kappa
\]
If $[\lambda_{l}]_{0}=0$, by the Eq.\ref{4.3.4}, then
\[
[\lambda_{l}]_{1}^{2}+[\lambda_{l}]_{2}^{2}+[\lambda_{l}]_{3}^{2}-2[\lambda_{l}]_{3}+2
=[\lambda_{l}]_{1}^{2}+[\lambda_{l}]_{2}^{2}+([\lambda_{l}]_{3}-1)^{2}+1=0
\]
Obviously, a real solution of this equation does not exist.

It is easy to see $\kappa$ is a left eigenvalue of $\mathbf{A}$. In conclusion, $\lambda_{l}=\pm1+\kappa$ or $\lambda_{l}=\kappa$.

\begin{Exa}\label{4.4}
Let
\[
\mathbf{A}=
\left[
\begin{array}{ccc}
\hbar&0&\hbar
\\
0&\jmath&0
\\
\hbar&0&\kappa
\end{array}
\right]
\]
\end{Exa}
Then
\[
\mathcal{P}_{1}(\mathbf{A}-\lambda_{l} \mathbf{I})=
\]
\begin{align}
\begin{bmatrix}\begin{smallmatrix}
-[\lambda_{l}]_{0}&-1+[\lambda_{l}]_{1}&[\lambda_{l}]_{2}&[\lambda_{l}]_{3}&0&0&0&0&0&-1&0&0
\\
1-[\lambda_{l}]_{1}&-[\lambda_{l}]_{0}&[\lambda_{l}]_{3}&-[\lambda_{l}]_{2}&0&0&0&0&1&0&0&0
\\
-[\lambda_{l}]_{2}&-[\lambda_{l}]_{3}&-[\lambda_{l}]_{0}&-1+[\lambda_{l}]_{1}&0&0&0&0&0&0&0&-1
\\
-[\lambda_{l}]_{3}&[\lambda_{l}]_{2}&1-[\lambda_{l}]_{1}&-[\lambda_{l}]_{0}&0&0&0&0&0&0&1&0
\\
0&0&0&0&-[\lambda_{l}]_{0}&[\lambda_{l}]_{1}&-1+[\lambda_{l}]_{2}&[\lambda_{l}]_{3}&0&0&0&0
\\
0&0&0&0&-[\lambda_{l}]_{1}&-[\lambda_{l}]_{0}&[\lambda_{l}]_{3}&1-[\lambda_{l}]_{2}&0&0&0&0
\\
0&0&0&0&1-[\lambda_{l}]_{2}&-[\lambda_{l}]_{3}&-[\lambda_{l}]_{0}&[\lambda_{l}]_{1}&0&0&0&0
\\
0&0&0&0&-[\lambda_{l}]_{3}&-1+[\lambda_{l}]_{2}&-[\lambda_{l}]_{1}&-[\lambda_{l}]_{0}&0&0&0&0
\\
0&-1&0&0&0&0&0&0&-[\lambda_{l}]_{0}&[\lambda_{l}]_{1}&[\lambda_{l}]_{2}&-1+[\lambda_{l}]_{3}
\\
1&0&0&0&0&0&0&0&-[\lambda_{l}]_{1}&-[\lambda_{l}]_{0}&-1+[\lambda_{l}]_{3}&-[\lambda_{l}]_{2}
\\
0&0&0&-1&0&0&0&0&-[\lambda_{l}]_{2}&1-[\lambda_{l}]_{3}&-[\lambda_{l}]_{0}&[\lambda_{l}]_{1}
\\
0&0&1&0&0&0&0&0&1-[\lambda_{l}]_{3}&[\lambda_{l}]_{2}&-[\lambda_{l}]_{1}&-[\lambda_{l}]_{0}
\end{smallmatrix}\end{bmatrix}
\end{align}
Let $\mathbf{B}=\mathcal{P}_{1}(\mathbf{A}-\lambda_{l} \mathbf{I})$, From the Theorem \ref{2.7}, let $\lambda_{l}\neq \jmath$, we can take
\begin{equation}\label{4.4.1}
\det\mathbf{B}(1:9,4:12)=-([\lambda_{l}]_{1}+[\lambda_{l}]_{3}-2[\lambda_{l}]_{1}[\lambda_{l}]_{3}-1)\beta=0
\end{equation}
\begin{equation}\label{4.4.2}
\det\mathbf{B}(1:8,10,4:12)=-([\lambda_{l}]_{0}-[\lambda_{l}]_{2}+2[\lambda_{l}]_{1}[\lambda_{l}]_{2})\beta=0
\end{equation}
\begin{equation}\label{4.4.3}
\det\mathbf{B}(1:8,11,4:12)=-([\lambda_{l}]_{0}^{2}-[\lambda_{l}]_{1}^{2}+[\lambda_{l}]_{1}+[\lambda_{l}]_{2}^{2}+[\lambda_{l}]_{3}^{2}-[\lambda_{l}]_{3}+1)\beta=0
\end{equation}
\begin{equation}\label{4.4.4}
\det\mathbf{B}(1:8,12,4:12)=([\lambda_{l}]_{0}+[\lambda_{l}]_{2}-2[\lambda_{l}]_{0}[\lambda_{l}]_{1})\beta=0
\end{equation}
and
\[
\beta=[\lambda_{l}]_{0}^{2}+[\lambda_{l}]_{1}^{2}+([\lambda_{l}]_{2}-1)^{2}+[\lambda_{l}]_{3}^{2}.
\]
By $\lambda_{l}\neq \jmath$, we can get\\
\[
\beta\neq 0
\]
then\\
\begin{equation}\label{4.4.4.1}
[\lambda_{l}]_{1}+[\lambda_{l}]_{3}-2[\lambda_{l}]_{1}[\lambda_{l}]_{3}-1=0
\end{equation}
\begin{equation}\label{4.4.4.2}
[\lambda_{l}]_{0}-[\lambda_{l}]_{2}+2[\lambda_{l}]_{1}[\lambda_{l}]_{2}=0
\end{equation}
\begin{equation}\label{4.4.4.3}
[\lambda_{l}]_{0}^{2}-[\lambda_{l}]_{1}^{2}+[\lambda_{l}]_{1}+[\lambda_{l}]_{2}^{2}+[\lambda_{l}]_{3}^{2}-[\lambda_{l}]_{3}+1=0
\end{equation}
\begin{equation}\label{4.4.4.4}
[\lambda_{l}]_{0}+[\lambda_{l}]_{2}-2[\lambda_{l}]_{0}[\lambda_{l}]_{1}=0
\end{equation}
From Eq.\ref{4.4.4.2} and \ref{4.4.4.4}, we get
\[
[\lambda_{l}]_{2}(2[\lambda_{l}]_{1}^{2}-2[\lambda_{l}]_{1}+1)=0
\]
i.e.
\[
[\lambda_{l}]_{2}=0
\]
And by the Eq.\ref{4.4.4.4}, get $[\lambda_{l}]_{0}=0$. Substituting $[\lambda_{l}]_{0}=0$ and $[\lambda_{l}]_{2}=0$ into Eq.\ref{4.4.4.3}, we get
\begin{equation}\label{4.4.4.5}
-[\lambda_{l}]_{1}^{2}+[\lambda_{l}]_{1}+[\lambda_{l}]_{3}^{2}-[\lambda_{l}]_{3}+1=0
\end{equation}

If $[\lambda_{l}]_{3}=\frac{1}{2}$, by the Eq.\ref{4.4.4.1}, we get $\frac{1}{2}=0$, this is contradictory.

If $[\lambda_{l}]_{3}\neq\frac{1}{2}$, then
\begin{equation}\label{4.4.4.6}
[\lambda_{l}]_{1}=\frac{1-[\lambda_{l}]_{3}}{1-2[\lambda_{l}]_{3}}
\end{equation}
Substituting the Eq.\ref{4.4.4.6} into the Eq.\ref{4.4.4.5}, then
\[
-(\frac{1-[\lambda_{l}]_{3}}{1-2[\lambda_{l}]_{3}})^{2}+\frac{1-[\lambda_{l}]_{3}}{1-2[\lambda_{l}]_{3}}+[\lambda_{l}]_{3}^{2}-[\lambda_{l}]_{3}+1=0
\]
i.e.
\begin{equation}\label{4.4.4.7}
[\lambda_{l}]_{3}^{4}-8[\lambda_{l}]_{3}^{3}+10[\lambda_{l}]_{3}^{2}-6[\lambda_{l}]_{3}+1=0
\end{equation}
Solving the Eq.\ref{4.4.4.7} (Polynomial), we get two real roots
\[
[\lambda_{l}]_{3}=\frac{1}{2}(1\pm\sqrt{-2+\sqrt{5}}\,)
\]
and two complex roots is omitted
\[
[\lambda_{l}]_{3}=\frac{1}{2}(1\pm \hbar\sqrt{-2+\sqrt{5}}\,)
\]
 by the Eq.\ref{4.4.4.6}, get
\[
[\lambda_{l}]_{1}=\frac{\frac{1}{2}\mp\sqrt{-2+\sqrt{5}}}{\mp\sqrt{-2+\sqrt{5}}},\,\,\,\,when\,\,\,\, [\lambda_{l}]_{3}=\frac{1}{2}(1\pm\sqrt{-2+\sqrt{5}}\,)
\]
It is easy to see $\jmath$ is a left eigenvalue of $\mathbf{A}$. In conclusion,
\[
\lambda_{l}=\jmath,\,\,\,\,\lambda_{l}=\frac{\frac{1}{2}\mp\sqrt{-2+\sqrt{5}}}{\mp\sqrt{-2+\sqrt{5}}}\hbar+\frac{1}{2}(1\pm\sqrt{-2+\sqrt{5}}\,)\kappa
\]
\begin{Exa}
Let
\[
\mathbf{A}=
\left[
\begin{array}{ccc}
a_{11}&0&0
\\
0&a_{22}&a_{23}
\\
0&a_{32}&a_{33}
\end{array}
\right].
\]
\end{Exa}
It is not hard to see that $a_{ij}$ and the left eigenvalues of $\left[\begin{array}{cc}
a_{22}&a_{23}
\\
a_{32}&a_{33}
\end{array}\right]$ are all the left eigenvalues of $\mathbf{A}$.

Since a $2\times2$ quaternion matrix may have one, two or an infinite number of
left eigenvalues \cite{huang2001left}, so a $3\times3$ quaternion matrix can have one, two, three or an infinite number of left eigenvalues. Similarly, a $n\times n$ quaternion matrix may have $1,2,\cdots,n-1,n$ or an infinite number of left eigenvalues. In addition, by Theorem \ref{2.7} $(2)$, by analysing the number of roots of specified four polynomial equations, obtains the number of left eigenvalues.

\begin{Rem}\label{4.4.4.4.4}
Numerical algorithms can be used to solve the specified four polynomial equations (Theorem \ref{2.7} $(2)$) to find the left eigenvalues of $m$ by $m$ quaternion matrix.
\end{Rem}

\section{Conclusion and Discussion}\label{S5}
In this paper, we introduce a method to compute the left eigenvalues of quaternion matrix based on the matrix representation of quaternion. We obtain four real polynomial equations with four variables which are equivalent to the generalized characteristic polynomial, the left eigenvalues could be found via solving special polynomial equations. In addition, while $\mathbf{A}\in\mathbb{H}^{m\times m}$ may have infinite number of left eigenvalues, the norm of left eigenvalues $\lambda_{l} $ of $\mathbf{A}$ is dominated by the norm of the right eigenvalues $\lambda_{r}$ of $\mathbf{A}$ .


Further theoretical analysis of the special polynomial equations about the left eigenvalues and finding much more potential applications of matrix representation will be left for future work.


\bibliography{mybibfile}

\end{document}